# Using Hermite Function for Solving Thomas-Fermi Equation

F. Bayatbabolghani, K. Parand

*Abstract*—In this paper, we propose Hermite collocation method for solving Thomas-Fermi equation that is nonlinear ordinary differential equation on semi-infinite interval. This method reduces the solution of this problem to the solution of a system of algebraic equations. We also present the comparison of this work with solution of other methods that shows the present solution is more accurate and faster convergence in this problem.

## I. Introduction

### A. Spectral Method

MANY of the current science and engineering problems are set in unbounded domains. In the context of spectral methods such as collocation and Galerkin methods [1], a number of approaches for treating unbounded domains have been proposed and investigated. The most common method is the use of polynomials that are orthogonal over unbounded domains, such as the transformed Hermite and Laguerre spectral method [2]–[9].

Guo [10]-[13] proposed a method that proceeds by mapping the original problem in an unbounded domain to a problem in a bounded domain, and then using suitable Jacobi polynomials such as Gegenbauer polynomials to approximate the resulting problems. The Jacobi polynomials are a class of classical orthogonal polynomials and the Gegenbauer polynomials, and thus also the Legendre and Chebyshev polynomials, are special cases of these polynomials which have been used in sevral literatures for solving some problems [14], [15].

On more approach is replacing infinite domain with $[-L, L]$ and semi-infinite interval with $[0, L]$ by choosing $L$, sufficiently large. This method is named domain truncation [16].

There is another effective direct approach for solving such problems which is based on rational approximations. Christov [17] and Boyd [18], [19] developed some spectral methods on unbounded intervals by using mutually orthogonal systems of rational functions. Boyd [18] defined a new spectral basis, named rational Chebyshev functions on the semi-infiniteinterval, by mapping to the Chebyshev polynomials. Guo et al. [20] introduced a new set of rational Legendre functions which are mutually orthogonal in $L^2(0, \infty)$. They applied a spectral scheme using the rational Legendre functions for solving the Korteweg-de Vries equation on the half-line. Boyd et al. [21] applied pseudospectral methods on a semi-infinite interval and compared rational Chebyshev, Laguerre and mapped Fourier sine methods.

Parand et al. [22]–[27], applied spectral method to solve nonlinear ordinary differential equations on semi-infinite intervals. Their approach was based on rational tau and collocation methods.

### B. Introduction of the Problem

Thomas-Fermi equation is one of the most important nonlinear ordinary differential equations that occurs in semi-infinite interval, as following [25], [28], [29]:

$$\frac{d^2 y}{dx^2} = \frac{1}{\sqrt{x}} y^{\frac{3}{2}}(x), \quad (1)$$

which appears in the problem of determining the effective nuclear charge in heavy atoms. Also, it has following Boundary conditions:

$$y(0) = 1, \quad y(\infty) = 0. \quad (2)$$

The Thomas-Fermi equation is useful for calculating formfactors and for obtaining effective potentials which can be used as initial trial potentials in self-consistent field calculations. The problem has been solved by different techniques [25], [30]–[39].

References [31]–[33] used perturbative approach to determine analytic solutions for Thomas-Fermi equation. Adomian [34] applied the decomposition method for solving Thomas-Fermi equation and then Wazwaz [35] proposed a nonperturbative approximate solution to this equation by using the modified decomposition method in a direct manner without any need to a perturbative expansion or restrictive assumptions. Liao [36] solved Thomas-Fermi equation by homotopy analysis method. Khan [37], used the homotopy analysis method (HAM) with a new and better transformation which improved the results in comparison with Liao's work. In [38], the quasilinearization approach was applied for solving (1). This method approximated the solution of a nonlinear differential equation by treating the nonlinear terms as a perturbation about the linear ones, and unlike perturbation theories is not based on the existence of some kind of a small parameter. Ramos [39] presented two piecewise quasilinearization methods for the numerical solution of (1).

F. Bayatbabolghani is with Department of Computer Science, shahid Beheshti University, Tehran, 1983963113Iran (Phone: +98 938 4115371; e-mail: fattaneh.bayat@gmail.com).

K. Parand is with Department of Computer Science, Shahid Beheshti University, Tehran,1983963113 Iran (e-mail: k_parand@sbu.ac.ir).

Both methods were based on the piecewise linearization of ordinary differential equations [25]. In addition, Parand [25] Solved Thomas-Fermi equation by Rational Chebyshev pseudospectral approach.

In this paper, we are going to solve Thomas-Fermi equation numerically by using the transformed Hermite functions via collocation method. we also have a comparison with a numerical solution.

Sections II reviews the desirable properties of Hermit functions with solution of the problem with collocation method by these functions, respectively. In Section IV we describe our results via tables and figures. Finally, concluding remarks will be presented in Section V.

## II. HERMITE FUNCTIONS

This section are devoted to elaborate the properties of Hermite functions. First of all, we should mention Hermite polynomials are generally not suitable in practice due to their wild asymptotic behavior at infinities [40]; therefore, we shall consider the Hermite functions. The normalized Hermite functions of degree $n$ is defined by [41]

$$\tilde{H}_n = \frac{1}{\sqrt{2^n n!}} e^{\frac{-x^2}{2}} H_n(x), \quad n \geq 0, x \in \mathbb{R}. \quad (3)$$

That $\{\tilde{H}_n\}$ is an orthogonal system in $L^2(\mathbb{R})$.

In the contrary to Hermite Polynomials, the Hermite functions are well behaved with the decay property:

$$|\tilde{H}_n(x)| \to 0, \quad as \ |x| \to \infty, \quad (4)$$

and, the three-term recurrence relation of Hermite functions implies [41]

$$\tilde{H}_{n+1}(x) = x\sqrt{\frac{2}{n+1}}\tilde{H}_n(x) - \sqrt{\frac{n}{n+1}}\tilde{H}_{n-1}(x), \quad n \geq 1, \quad (5)$$

$$\tilde{H}_0(x) = e^{\frac{-x^2}{2}},$$

$$\tilde{H}_1(x) = \sqrt{2} x e^{\frac{-x^2}{2}}.$$

For more details you can study [41]–[43].

Steady flow problem is defined on the interval $(0, +\infty)$, but Hermite functions are defined on the interval $(-\infty, +\infty)$. One of the approaches to construct an approximation on the interval $(0, +\infty)$ is using mapping that is changing variable of the form [41]

$$w = \Phi(z) = \frac{1}{k}\ln(z), \quad (6)$$

where $k$ is a constant.

he transformed Hermite functions are

$$\hat{H}_n(x) \equiv \tilde{H}_n(x) o \Phi(x) = \tilde{H}_n(\Phi(x)), \quad (7)$$

The inverse map of $w = \Phi(z)$ is

$$z = \Phi^{-1}(w) = e^{kw}. \quad (8)$$

Therefore, we may define the inverse images of the spaced nodes $\{x_j\}_{x_j=-\infty}^{x_j=+\infty}$ as [41]

$$\Gamma = \{\Phi^{-1}(t) : -\infty < t < +\infty\} = (0, +\infty), \quad (9)$$

and

$$\tilde{x}_j = \Phi^{-1}(x_j) = e^{x_j}, \quad j = 0,1,2,\ldots \quad (10)$$

Let $w(x)$ denotes a non-negative, integrable, real-valued function over the interval $\Gamma$, We define [41]

$$L^2_w(\Gamma) = \{v : \Gamma \to \mathbb{R} \mid v \text{ is measurable and } \|v\|_w < \infty\}, \quad (11)$$

where

$$\|v\|_w = \left(\int_0^\infty |v(x)|^2 w(x) dx\right)^{\frac{1}{2}}, \quad (12)$$

is the norm induced by the inner product of the space $L^2_w(\Gamma)$ [41],

$$<u,v>_w = \int_0^\infty u(x)v(x)w(x)dx. \quad (13)$$

Thus, $\{\hat{H}_n(x)\}_{n \in N}$ denotes a system which is mutually orthogonal

$$\langle \hat{H}_n(x), \hat{H}_m(x) \rangle_{w(x)} = \sqrt{\pi}\delta_{nm}. \quad (14)$$

This system is complete in $L^2_w(\Gamma)$. Therefore, for any function $f \in L^2_w(\Gamma)$ the following expansion holds [41]

$$f(x) \cong \sum_{k=-N}^{+N} f_k \hat{H}_k(x), \quad (15)$$

with

$$f_k = \frac{\langle f(x), \hat{H}_k(x) \rangle_{w(x)}}{\left\|\hat{H}_k(x)\right\|^2_{w(x)}} \quad (16)$$

Now we define an orthogonal projection based on the transformed Hermite function as given below [41].

Let
$$\hat{H}_N = span\{\hat{H}_0(x), \hat{H}_1(x), \ldots, \hat{H}_n(x)\}. \quad (17)$$

The $L^2(\Gamma)$-orthogonal projection $\hat{\xi}_N : L^2(\Gamma) \to \hat{H}_N$ is a mapping in a way that for any $y \in L^2(\Gamma)$ [41],

$$\langle \hat{\xi}_N y - y, \phi \rangle = 0 \quad \forall \phi \in \hat{H}_N, \quad (18)$$

or equivalently,

$$\hat{\xi}_N y(x) = \sum_{i=0}^{N} \hat{a}_i \hat{H}_i(x). \quad (19)$$

### III. Solving the Problem with Hermite Functions

For solving Thomas-Fermi, we used $\frac{1}{k}\ln(x)$ for changing variable. Also, because of boundary conditions, we set following function:

$$p(x) = \frac{1}{1 + \lambda x + x^2}, \quad (20)$$

and $\lambda$ is constant.

Finally, we have

$$\hat{\xi}_N f(x) = p(x) + \hat{\xi}_N f(x). \quad (21)$$

that

$$\hat{\xi}_N f(x) = \sum_{i=0}^{N} \hat{a}_i \hat{H}_i(x). \quad (22)$$

To find the unknown coefficients $\hat{a}_i$'s, we substitute the truncated series $\hat{\xi}_N f(x)$ into (1). Also, we define Residual function of the form

$$Res(x) = p''(x) + \hat{\xi}_N f''(x) - x^{-\frac{1}{2}} \{p(x) + \hat{\xi}_N f(x)\}^{\frac{3}{2}} = 0. \quad (23)$$

By applying $x$ in (23) with the $N$ collocation points which are roots of transformed Hermite function, we have $N$ equations that generates a set of $N$ nonlinear equations. Now, all of these equations can be solved by Newton method for the unknown coefficients.

### IV. Result

The initial slope $y'(0)$ of the Thomas-Fermi equation is calculated by Kobayashi [44] as $y'(0) = -1.588071$. Table I shows the approximations of $y(x)$ and $y'(0)$ obtained by the present method for $N = 15$, $k = 0.9$ and $\lambda = 1.588071$, and those obtained by Liao [45] and Kobayashi [44].

TABLE I
COMPARISON BETWEEN TRANSFORMED HERMITE FUNCTION AND LIAO [45] WITH $N = 15$, $k = 0.9$ AND $\lambda = 1.588071$.

| $x$ | Present method | Liao [45] |
|---|---|---|
| 00.25 | 0.754795330 | 0.755202000 |
| 00.50 | 0.606658908 | 0.606987000 |
| 00.75 | 0.502110510 | 0.502347000 |
| 01.00 | 0.423811203 | 0.424008000 |
| 01.25 | 0.363027725 | 0.363202000 |
| 02.00 | 0.242918233 | 0.243009000 |
| 02.25 | 0.215819818 | 0.215895000 |
| 02.50 | 0.192917948 | 0.192984000 |
| 02.75 | 0.173379623 | 0.173441000 |
| 03.00 | 0.156573773 | 0.156633000 |
| 03.25 | 0.142013368 | 0.142070000 |
| 03.50 | 0.129316613 | 0.129370000 |
| 03.75 | 0.118180209 | 0.118229000 |
| 04.00 | 0.108360441 | 0.108404000 |
| 08.00 | 0.036580427 | 0.036587300 |
| 15.00 | 0.010803774 | 0.010805400 |
| 20.00 | 0.005792831 | 0.005784940 |
| 30.00 | 0.002252634 | 0.005784940 |
| $y'$ | Present method | Kobayashi [44] |
| $y'(0)$ | -1.588071 | -1.588071 |

Fig. 1 shows the resulting graph of Thomas-Fermi for $N = 15$, $k = 0.9$ and $\lambda = 1.588071$ which tends to zero as x increases by boundary condition $y(\infty) = 0$. It is compared with Liao's results that show by square.

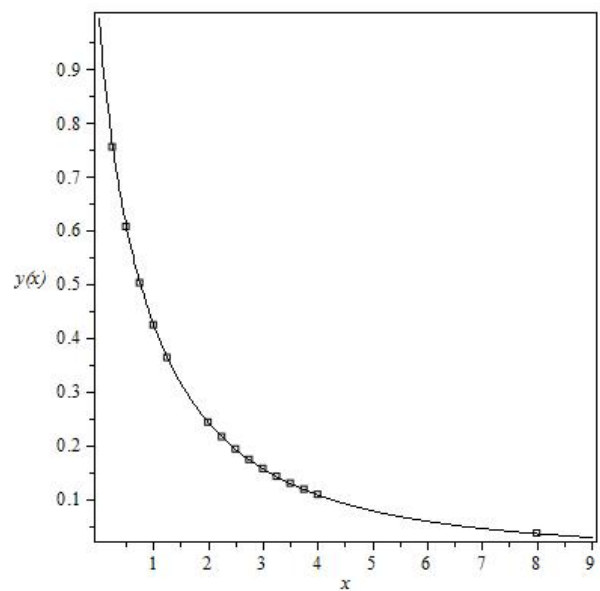

Fig. 1 Graph of comparison numerical approximate $y(x)$ by transformed Hermite functions to Liao [45] with $N = 15$, $k = 0.9$ and $\lambda = 1.588071$

## V. Conclusions

In the above discussion, we applied the collocation method to solve the Thomas-fermi equation that is defined in a semi-infinite interval which has singularity at $x = 0$ and its boundary condition occurs in infinity. Our scheme was based on transformed Hermite function that solved the non-linear differential equations on the semi-infinite domain without truncating it to a finite domain. Transformed Hermite function was proposed to provide simple way to improve the convergence of the solution by collocation method. Finally, we reported our numerical finding and demonstrated the present solution was highly accurate.